\documentclass[10pt,epsfig]{article}   
\usepackage{amsbsy,amssymb,amscd,amsfonts,latexsym,amstext,delarray,   
amsmath,epsfig}  \headsep=15pt   
  
\input xypic   
   
\newtheorem{thm}{Th\'eor\`eme}[section]

\newtheorem{defn}[thm]{D\'efinition}

\numberwithin{equation}{section}

\def\R{{\mathbb R}}

\def\H{{\mathbb H}}

\def\mX{{\mathfrak X}}

\def\C{{\mathbb C}}   
\def\Z{{\mathbb Z}}   
\def\N{{\mathbb N}}   
\def\HH{{\mathbb H}}   
\def\P{{\mathbb P}}   
   
\def\PSL{{\rm PSL}}  
   
\def\SL{{\rm SL}}   
\def\Aut{{\rm Aut}}  
\def\Tr{{\rm Tr}}  
\def\w{{\text w}}   
   
\def\XX{{{X^*}}}

\newcommand{\cf}{{\it voir}\ }

\begin{document}   
\title{Triplets spectraux en g\'eom\'etrie d'Arakelov   
}   
   
\author{Caterina Consani\thanks{Partiellement support\'ee par la
bourse du NSERC num\'ero 72016789}\\ D\'epartement de   
Math\'ematiques \\   
Universit\'e de Toronto, Canada \\ kc\@@math.toronto.edu 
\and Matilde Marcolli   
\thanks{Partiellement support\'ee par la bourse Sofja Kovalevskaja de
l'Humboldt Foundation} \\    
Max--Planck Institut f\"ur Mathematik  \\   
Bonn, Allemagne \\ marcolli\@@mpim-bonn.mpg.de}   
\date{}  
  
\maketitle  

\noindent {\bf R\'esum\'e} Dans cette note nous employons la th\'eorie
des triplets spectraux de Connes pour rapprocher le mod\`ele de Manin
du graphe dual de la fibre \`a l'infini d'une surface d'Arakelov et la
cohomologie du c\^one de la monodromie locale.

\bigskip

\section*{Spectral triples in Arakelov geometry}

\noindent {\bf Abstract} 
In this note, we use Connes' theory of spectral triples to provide a
connection between Manin's model of the dual
graph of the fiber at infinity of an Arakelov surface and the
cohomology of the mapping cone of 
the local monodromy.

\bigskip

\subsection*{Abridged English version}

 In Arakelov theory a completion of an arithmetic
surface is achieved by enlarging the group of divisors by formal
linear combinations of the ``closed fibers at infinity''.
Manin in \cite{Man} described the dual graph of any such closed fiber
in terms of an infinite tangle of bounded geodesics in a hyperbolic
handlebody $\mX_\Gamma = \Gamma \backslash \H^3$, uniformized by a
Schottky group $\Gamma \subset \PSL(2,\C)$. In this note we consider
arithmetic surfaces over the ring of integers in a number field, with
fibers of genus $g\geq 2$. We use Connes' theory of spectral triples 
to relate the hyperbolic geometry of the handlebody to the cohomology of
the cone of the local monodromy $N$ at arithmetic infinity as
introduced in \cite{KC}. 
First, we construct a spectral triple, where the non--commutative space is
given by the reduced ${\rm C}^*$--algebra of the Schottky group acting
on the cohomology of the cone via a representation induced by the
presence of a polarized Lefschetz module structure. 
In this setting we recover the alternating product of the archimedean 
factors from a zeta function of the spectral triple. 
Then, we introduce a second spectral triple, which is related
to Manin's description of the dual graph of the fiber at infinity.
Here the non--commutative space is a ${\rm C}^*$--algebra representing
the ``reduction mod infinity'' and acting on a ``dynamical
homology and cohomology'' pair, defined in terms of the bounded
geodesics in the handlebody and of a dynamical system $T$. 
The operator $\Phi$, that represents the ``logarithm of a Frobenius--type
operator'' on the archimedean cohomology of \cite{Den}, gives the
Dirac operator on these spectral triples. We show that the archimedean
cohomology 
embeds in the dynamical cohomology, compatibly with the action of a
real Frobenius $\bar F_\infty$, so that the duality isomorphism on the
cohomology of the cone of $N$ corresponds to the pairing of dynamical
homology and cohomology. 

A detailed version of the results presented
in this note is contained in \cite{CM}.

\bigskip

\section{Introduction} 
  
En th\'eorie d'Arakelov une compl\'etion d'une surface arithm\'etique  
est r\'ealis\'ee par l'\'elargissement du groupe des diviseurs avec une
combinaison lin\'eaire formelle    
des ``fibres ferm\'ees \`a l'infini''.  
Dans \cite{Man}, Manin a d\'ecrit le graphe dual d'une telle fibre ferm\'ee  
en employant un enlacement infini de g\'eod\'esiques limit\'ees dans
une vari\'et\'e hyperbolique   
$\mX_\Gamma = \Gamma \backslash \H^3$, uniformis\'ee par un  
groupe de Schottky $\Gamma \subset \PSL(2,\C)$. Dans cette note nous
consid\'erons   
des surfaces arithm\'etiques d\'efinies sur l'anneau des entiers d'un
corps de nombres, avec   
des fibres de genre $g\geq 2$. Nous employons la th\'eorie des triplets
spectraux de Connes     
pour rapprocher la g\'eom\'etrie hyperbolique de la vari\'et\'e
avec la cohomologie du   
c\^one de la monodromie locale $N$ \`a l'infini aritm\'etique qui a
\'et\'e introduite dans \cite{KC}.    
D'abord, nous construisons un triplet spectral, o\`u l'espace
non--commutatif est donn\'e   
par la ${\rm C}^*$--alg\`ebre r\'eduite du groupe de Schottky qui agit  
sur la cohomologie du c\^one \`a travers la repr\'esentation d\'efinie par la  
pr\'esence d'une structure de Lefschetz module polaris\'ee.   
Dans ce cadre nous retrouvons le produit altern\'e des facteurs
archim\'ediens   
d'une fonction zeta du triplet spectral.   
Puis, nous introduisons un deuxi\`eme triplet spectral, qui est li\'e  
\`a la description de Manin du graphe dual de la fibre \`a l'infini.  
Ici l'espace non--commutatif est une ${\rm C}^*$--alg\`ebre qui repr\'esente  
la ``r\'eduction modulo infini'' et qui agit sur un accouplement d'
``homologie et cohomologie dynamique'', d\'efinie   
avec les  g\'eod\'esiques limit\'ees dans la vari\'et\'e hyperbolique
et \`a travers l'usage d'un syst\`eme dynamique $T$.    
L'op\'erateur $\Phi$, qui repr\'esente le ``logarithme de Frobenius''   
sur la cohomologie archim\'edienne de \cite{Den}, donne l'op\'erateur de  
Dirac sur ces triplet spectraux. Nous prouvons que la cohomologie
archim\'edienne est plong\'ee
dans la cohomologie dynamique, de mani\`ere compatible avec l'action de  
Frobenius r\'eelle $\bar F_\infty$, de sorte que l'isomorphisme de
dualit\'e sur la    
cohomologie du c\^one de $N$ correspond \`a l'accouplement
d'homologie et cohomologie dynamique.   
  
Une version d\'etaill\'ee des r\'esultats pr\'esent\'es dans cette  
note est contenue dans \cite{CM}. \vspace{.2in}  

\section{R\'esultats}
  
Soit $X_{/\kappa}$ une courbe projective et lisse d\'efinie sur $\kappa =  
\C$ or $\R$. Pour $a,b\in\N$, nous notons    
$(A^{a,b}\oplus A^{b,a})_\R$ le groupe ab\'elien des formes
diff\'erentielles r\'eelles    
(analytiques ou ${\rm C}^\infty$) sur $X_{/\kappa}$ de type   
$(a,b) + (b,a)$.     
   
Pour $p\in\Z$, l'expression $(A^{a,b}\oplus A^{b,a})_\R(p)$ signifie   
le $p$-\`eme Hodge-Tate twist de $(A^{a,b}\oplus A^{b,a})_\R$.   
   
Soient $i,j,k\in\Z$; nous consid\'erons le complexe suivant  
\begin{equation}\label{complex}   
K^{i,j,k} =   
\begin{cases}   
{\displaystyle \bigoplus_{\begin{subarray}{l} a+b=j+1\\|a-b|\le   
2k-i\end{subarray}} (A^{a,b}\oplus A^{b,a})_\R(\frac{1+j-i}{2})}   
&\text{si $1+j-i \equiv~0(2),~k\ge\text{max}(0,i)$} \\   
0 &\text{sinon.}   
\end{cases}   
\end{equation}   
   
Sur $K^{i,j,k}$ on d\'efinit les diff\'erentielles    
$ d': K^{i,j,k} \to K^{i+1,j+1,k+1}$ et $d'': K^{i,j,k} \to   
K^{i+1,j+1,k}$, avec $d'=  
\partial+\overline\partial$ et $d''=   
\sqrt{-1}(\overline\partial-\partial)$.   
Ceux-ci  satisfont  ${d'}^2 = 0 = {d''}^2$ (\cf\cite{KC}  
Lemma~4.2).   
Nous consid\'erons aussi les morphismes  
\[ N: K^{i,j,k} \to K^{i+2,j,k+1},\quad N(f) =  
(2\pi\sqrt{-1})^{-1}f \  \text{ et }   \ l: K^{i,j,k} \to  
K^{i,j+2,k},\quad l(f) = (2\pi\sqrt{-1})f\wedge\omega, \]  
o\`u $N$ peut \^etre regard\'e comme le logarithme de la  {\it
monodromie locale \`a l'infini},   
 et $l$ est l'{\it homomorphisme de Lefschetz}, avec  
$\omega$ la $(1,1)$-forme r\'eelle fondamentale (ferm\'ee) sur
$X_{/\kappa}$.    
Ces endomorphismes commutent avec $d'$ et $d''$ et   
satisfont $[l,N] = 0$ (\cf\cite{KC}). On d\'efinit $K^{i,j} = \oplus_k  
K^{i,j,k}$ et on \'ecrit     
$K^\ast = \oplus_{i+j=\ast}K^{i,j}$ pour indiquer le complexe simple   
dot\'e de la diff\'erentielle total $d = d'+d''$ et avec l'action   
des op\'erateurs $N$ et $l$.  
  
Le complexe diff\'erentiel $K^{\cdot,\cdot}$  
avec les op\'erateurs $N$ et $l$ est un {\em module de Lefschetz
bigradu\'e et polaris\'e}, avec    
la polarisation $\psi : K^{-i,-j,k}\otimes K^{i,j,k+i} \to \R(1)$   
d\'efinie par   
\[  
\psi(x,y) := \left( \frac{1}{2\pi \sqrt{-1}} \right)  
\epsilon(1-j) (-1)^k \int_{X_{/\kappa}(\C)} x \wedge Cy.   
\]  
Ici, pour $m\in\Z$: $\epsilon(m) := (-1)^{\frac{m(m+1)}{2}}$ et   
$C(x):= (\sqrt{-1})^{a-b}  x$ est l'op\'erateur de Weil pour $x$  
une forme diff\'erentielle de type $(a,b)$ (\cf\cite{Wells} \S V.1).  
La forme bilin\'eaire   
\begin{equation}\label{innerprod}  
\langle\cdot,\cdot\rangle: K\otimes K \to \R(1), \ \ \ \  
\langle x, y \rangle : = \psi\left( x, \sigma(\tilde w) y \right)  
\end{equation}  
est sym\'etrique et d\'efinie positive (\cf \cite{KC} Lemmas~4.2,~4.5,~4.6 et  
Proposition~4.7).  
La structure de module de Lefschetz  bigradu\'ee correspond \`a  
la repr\'esentation   
\begin{equation}\label{repSLSL}  
 \sigma: \SL(2,\R)\times \SL(2,\R) \to \Aut (K^{\cdot,\cdot})   
\end{equation}  
\[ \sigma\left\{ \left(\begin{array}{cc} a & 0 \\ 0 & a^{-1}  
\end{array} \right) , \left(\begin{array}{cc} b & 0 \\ 0 & b^{-1}  
\end{array} \right) \right\} (x) = a^i b^j x \ \ \ \text{ o\`u $x\in  
K^{i,j}$} \]  
\[ d\sigma \left\{ \left( \begin{array}{cc} 0 & 1 \\ 0 & 0 \end{array}
\right),   
0 \right\} = N \ \ \ \  d\sigma \left\{ 0, \left( \begin{array}{cc} 0  
& 1 \\ 0 & 0 \end{array}   
\right) \right\} = l. \]  
Nous consid\'erons le c\^one de l'application $N$  
\[  
{\rm Cone}(N)^{\cdot,\cdot} = {\rm Cone}(N: K^{\cdot,\cdot}\to K^{\cdot  
+2, \cdot}) := K^{\cdot,\cdot}[1] \oplus K^{\cdot +2,\cdot}, \ \ \  
D(a,b) = (-d(a),N(a)+d(b))   
\]  
et l'hyper-cohomologie $H^q(\XX) := \HH^q({\rm Cone}(N)^\cdot)$.   
Ces groupes ont une structure gradu\'e $H^q(\XX)= \oplus_{p\in \Z}  
gr_{2p}^\w H^q(\XX)$ de la forme  
\[ \begin{array}{l}  
H^0(\XX) = \oplus_{p\le 0}H^0(X_{/\C},\R(p))  \ \ \ \    
H^1(\XX) = \oplus_{p\le 0}H^1(X_{/\C},\R(p))\oplus   
\oplus_{p\ge 1}H^0(X_{/\C},\R(p-1))\\[2mm]  
  H^{2}(\XX) =\oplus_{p\le 1}H^2(X_{/\C},\R(p))\oplus \oplus_{p\ge  
2}H^1(X_{/\C},\R(p-1)) \ \ \ \ \        
H^3(\XX) = \oplus_{p\ge 2}H^2(X_{/\C},\R(p-1)). \end{array} \]  
Quand $\kappa = \R$ on obtient des r\'esultats similaires en prenant les
invariants  de   
la conjugaison de de Rham $\bar F_{\infty}$ (\cf \cite{CM} \S 2.)  
  
\medskip  
  
Soit $\Gamma\subset \PSL(2,\C)$ un groupe de Schottky de rang  
$g\geq 2$. Soit $\tilde \Gamma$ le $\Gamma$-stabilisateur   
de chacune de deux composantes connexes en $\P^1(\C)\setminus C$, o\`u  
$C$ est un quasi-cercle pour $\Gamma$ (\cf \cite{Bo}). Nous consid\'erons les  
deux groupes Fuchsiens de Schottky     
$G_i :=\{ \alpha_i \gamma \alpha_i^{-1}: \, \gamma \in \tilde \Gamma  
\}$, avec $\alpha_i$ les \'equivalences conformelles entre les deux  
composantes de $\P^1(\C)\setminus C$ et les deux h\'emisph\`eres en  
$\P^1(\C)\setminus \P^1(\R)$.   
Apr\`es un rel\`evement de $\tilde\Gamma$ de $\PSL(2,\C)$ \`a
$\SL(2,\C)$, nous consid\'erons   
la ${\rm C}^*$--alg\`ebre r\'eelle, r\'eduite ${\rm C}^*(\tilde\Gamma)$.

\begin{thm} \label{spectral3} Soit ${\mathcal X}$ une surface
arithm\'etique avec   
fibres de dimension un et de genre $g\ge 2$. Soit $\Gamma$ un
groupe de Schottky   
qui fixe l'uniformisation de la surface de Riemann correspondante
$X_{/\C}$ \`a une place archim\'edienne.   
Le produit \eqref{innerprod} induit un produit int\'erieur  
sur $H^\cdot (\XX)$. La repr\'esentation  
$\sigma_2 (\gamma) := \sigma\{ 1, \alpha_2 \gamma \alpha_2^{-1} \}$,   
obtenue par une restriction de \eqref{repSLSL} sur le
$\SL(2,\R)$-rel\`evement du   
groupe Fuchsien $G_2$, induit une repr\'esentation de ${\rm  
C}^*(\tilde\Gamma)$ sur la compl\'etion de Hilbert de $H^q(\XX)$ avec
le produit int\'erieur ci-dessus.  Ces donn\'ees,
avec l'op\'erateur   
\[ \Phi|_{gr_{2p}^\w H^q(\XX)} :=\left\{ \begin{array}{ll} p & q\geq 2p \\   
p-1 & q\leq 2p-1, \end{array} \right. \]   
d\'eterminent un triplet spectral 1--sommable $({\rm  
C}^*(\tilde\Gamma),H^\cdot (\XX),\Phi)$ \`a la Connes (\cf~\cite{Connes}).   
\end{thm}

Sur les sous-espaces   
\[ \begin{array}{l}  
H^-(\XX):= \oplus_{p\leq 0} gr_{2p}^\w H^0(\XX) \oplus \oplus_{p\le  
0} gr^\w_{2p} H^1(\XX) \oplus \oplus_{p\le 1} gr^\w_{2p} H^2(\XX),  
\\[3mm]   
 H^+(\XX):= \oplus_{p\ge 1} gr^\w_{2p} H^1(\XX) \oplus \oplus_{p\ge  
2} gr^\w_{2p} H^2(\XX) \oplus \oplus_{p\ge 2} gr^\w_{2p} H^3(\XX),   
\end{array} \]  
on d\'efinit des isomorphismes de dualit\'e $\delta=\oplus_{q=0}^2 \delta_q$  
qui sont produits par les puissances de la monodromie (\cf \cite{KC},  
Proposition~4.8). L'op\'erateur  
$$ \omega = \left(\begin{array}{cc} 0 & \delta^{-1} \\ \delta & 0  
\end{array} \right) $$  
satisfait $\omega^2= id$, $\omega^* =\omega$, $[\omega,a]=0$, pour chaque  
$a\in {\rm C}^*(\tilde\Gamma)$, et $(\Phi \omega + \omega  
\Phi)|_{H^q(\XX)} = q\cdot id$. Par cons\'equent, nous consid\'erons
dans la famille de fonctions zeta   
associ\'ees au triplet spectral  
\begin{equation}\label{zeta} \zeta_{a, P_- \Phi}(s,z):= \sum_{\lambda  
\in {\rm Spec}(P_- \Phi)} \Tr(a \Pi(\lambda, P_- \Phi)) (s-\lambda)^{-z},   
\end{equation}  
avec $a=\sigma_2(-id)$. Ici $P_-$ est la projection sur le sous-espace
$H^-(X^*)$ et $\Pi(\lambda, P_- \Phi)$ est la projection spectrale sur
le sous-espace propre de l'op\'erateur $P_- \Phi$ avec valeur 
propre $\lambda$.

\begin{thm} La fonction zeta \eqref{zeta} satisfait  
\[ \exp\left( -\frac{d}{dz} \zeta_{a,P_- \Phi/(2\pi)}(s/(2\pi),z)|_{z=0}  
\right)^{-1} = \frac{ L_\C (H^1(X_{/\C},\C),s) }{L_\C (H^0(X_{/\C},\C),s)  
\cdot L_\C (H^2(X_{/\C},\C),s) }, \]  
o\`u $L_\C (H^q(X_{/\C},\C),s)$ sont les facteurs $L$ archim\'ediens 
d\'efinis en \cite{Serre}.  
\end{thm}  
  
\medskip  
  
Pour un choix d'un ensemble de g\'en\'erateurs $\{ g_i \}_{i=1}^g$ du
groupe de   
Schottky $\Gamma$, nous consid\'erons l'ensemble ${\mathcal S}$ des {\em  
suites r\'eduites et doublement infinies}    
dans les $\{ g_i \}_{i=1}^{2g}$ ($g_{i+g}=g_i^{-1}$):  
\[ {\mathcal S}= \{\ldots a_{-m} \ldots a_{-1}  
a_0 a_1 \ldots a_\ell   
\ldots  \, \, | a_i \in \{ g_i \}_{i=1}^{2g}, \, \, a_{i+1} \neq a_i,  
\forall i\in \Z \}. \]  
L'{\em op\'erateur de d\'ecalage} $T$ agit sur ${\mathcal  
S}$ comme   
\[ T(\ldots a_{-m} \ldots a_{-1} a_0  
 a_1 \ldots  a_{\ell} \ldots ) = \ldots a_{-m+1}  \ldots  a_{0}   
a_1  a_2 \ldots  a_{\ell+1}  \ldots  \]  
Le couple $({\mathcal S},T)$ est un espace de Smale (\cf \cite{Putn}), et    
le tore de l'application $T$ est d\'efini comme   
${\mathcal S}_T := {\mathcal S} \times [0,1] / (x,0)\sim (Tx,1)$.   
  
Il y a une identification $H^1({\mathcal S}_T,\Z) \cong K_0({\rm  
C}({\mathcal S}) \rtimes_T \Z)$ sur la cohomologie de ${\mathcal S}_T$  
avec le $K_0$-groupe de la ${\rm C}^*$-alg\`ebre produit crois\'e pour  
l'action de $T$ sur ${\mathcal S}$ (\cf \cite{BoHa}).   
Ceci munit $H^1({\mathcal S}_T,\Z)$  
d'une filtration de groupes commutatifs libres $F_0  
\hookrightarrow F_1 \hookrightarrow \cdots F_n\hookrightarrow \cdots$,  
avec ${\rm rang}\, F_0=2g$ et  
${\rm rang}\, F_n = 2g(2g-1)^{n-1}(2g-2) +1$, pour $n\geq 1$, de
fa\c{c}on que   
$H^1({\mathcal S}_T,\Z)= \varinjlim_n F_n$ (\cf \cite{PaTu}).   
Le groupe d'homologie $H_1({\mathcal S}_T,\Z)$ a  
une filtration de groupes commutatifs libres ${\mathcal K}_N$, de sorte que   
$H_1({\mathcal S}_T,\Z) = \varinjlim_N {\mathcal K}_N$,    
avec ${\rm rang}( {\mathcal K}_N ) = (2g-1)^N +1$ pour $N$ pair et   
$(2g-1)^N + (2g-1)$ pour $N$ impair.  
  
\begin{defn}\begin{enumerate}\item
Nous d\'efinissons la {\em cohomologie dynamique} comme 
$$ H^1_{dyn} :=  \oplus_{p\leq 0} gr_{2p}^\Gamma H^1_{dyn}, \ \ \ \text{ 
o\`u } \ \ gr_{2p}^\Gamma  
H^1_{dyn} := {\rm Gr}_{-p} \otimes_{\R} \R(p), $$ 
pour ${\rm Gr}_n =(F_n/F_{n-1})\otimes_\Z \R$. Nous d\'efinissons le
sous-espace gradu\'e $ {\mathcal V}:= \oplus_{p\leq 0} gr_{2p}^\Gamma
{\mathcal V}$ de  $H^1_{dyn}$ avec $gr_{2p}^\Gamma {\mathcal V}$
engendr\'e par les   
\'el\'ements $(2\pi \sqrt{-1})^p \chi_{-p+1,k}$, pour $\chi_{n,k}:=[  
\chi_{{\mathcal S}^+(w_{n,k})} ]$ en ${\rm Gr}_{n-1}$.
Ici $\chi_{{\mathcal S}^+(w_{n,k})}$ est la fonction 
charact\'eristique du sous-espace topologique 
${\mathcal S}^+(w_{n,k})$ de 
$${\mathcal S}^+=\{ a_0 a_1 \ldots a_\ell \ldots \, | a_i \in 
\{ g_i \}_{i=1}^{2g}, \, a_{i+1} \neq a_i^{-1} \}$$
avec $a_0 \ldots a_n= w_{n,k} :=\underbrace{g_k g_k \ldots  
g_k}_{n-fois}$.  
  
\item Nous d\'efinissons l' {\em homologie dynamique} comme 
$$ H_1^{dyn}:=\oplus_{p\geq   
1} gr_{2p}^\Gamma H_1^{dyn}, \ \ \ \text{ o\`u } \ \ gr_{2p}^\Gamma
H_1^{dyn} :=  {\mathcal K}_{p-1}\otimes \R(p).$$ 
Nous d\'efinissons aussi ${\mathcal W}\subset   
H_1^{dyn}$ comme le sous-espace gradu\'e   
${\mathcal W}=\oplus_{p\geq 1} gr_{2p}^\Gamma {\mathcal W}$,  
o\`u $gr_{2p}^\Gamma {\mathcal W}$ est engendr\'e par les $2g$
\'el\'ements    
$(2\pi \sqrt{-1})^p \, \, \underbrace{g_k g_k \ldots  
g_k}_{p-fois}$. 
\end{enumerate}
\end{defn}

Il y a une involution (que nous notons encore   
$\bar F_\infty$) qui agit sur l'homologie et la cohomologie dynamiques
et qui est induit par le changement d'orientation.    
  
\begin{thm}\label{map-ar-dyn-2}  
Il y a des isomorphismes $\bar F_\infty$-\'equivariants $U$ et $\tilde U$,  
tels que le diagramme suivant soit commutatif (avec $p\le 0$):  
\begin{eqnarray*}  
\diagram  
gr_{2p}^\w H^1(\XX) \dto^U\rto^{\delta_1} &  gr_{2(-p+2)}^\w  
H^2(\XX)\dto^{\tilde U} \\  
gr_{2p}^\Gamma {\mathcal V}\rto^{{\mathcal D}} & gr_{2(-p+1)}^\Gamma   
{\mathcal W}.  
\enddiagram  
\end{eqnarray*}  
Ici $\delta_1$ est l'isomorphisme de dualit\'e de \cite{KC} qui est
induit par   
l'op\'erateur de monodromie $N$, et ${\mathcal D}$ est l'isomorphisme
de dualit\'e    
induit par l'accouplement d'homologie et cohomologie de ${\mathcal  
S}_T$. \end{thm}

L'alg\`ebre de Cuntz--Krieger ${\mathcal O}_A$ (\cf \cite{CuKrie}), avec  
la matrice \'el\'ementaire $A$ du sous-d\'eplacement de type fini $({\mathcal  
S},T)$, satisfait ${\mathcal O}_A\cong {\rm C}(\Lambda_\Gamma)\rtimes  
\Gamma$, o\`u $\Lambda_\Gamma\subset \P^1(\C)$ est l'ensemble limite du  
groupe de Schottky $\Gamma$.  
  
Les isom\'etries g\'en\'eratrices $S_i$ de ${\mathcal O}_A$ agissent sur
$F_n$ comme $(S_i   
h)(a_0\ldots a_n) = h(g_i^{-1} a_0\ldots a_{n-1}) \cdot \chi_{a_0\neq  
g_i}(a_0\ldots a_n)$. Cela d\'efinit une action de la ${\rm  
C}^*$-alg\`ebre (r\'eelle) ${\mathcal O}_A$ sur une compl\'etion de Hilbert
appropri\'ee de $H^1_{dyn}$. 
De plus, les fonctions $f\in {\rm C}_0(\H^3)$ d\'efinissent  
des op\'erateurs $(\tilde\rho(f) h) (a_0\ldots a_n) = f((a_0\ldots a_n)\cdot  
x_0) \cdot h(a_0\ldots a_n)$, avec $x_0 \in \H^3$ un point de base fix\'e.  
Ceci donne une repr\'esentation de  
$({\rm C}(\Lambda_\Gamma)\otimes {\rm C}_0(\H^3)) \rtimes  
\Gamma$ sur l'alg\`ebre des op\'erateurs limit\'es sur $H^1_{dyn}$.  
  
Les fonctions $f\in {\rm C}_b( \H^3,{\rm C}( \Lambda_\Gamma) )$
d\'efinissent des op\'erateurs   
$\tilde\rho(f) \, a_0\ldots a_p   
= f_{(a_0\ldots a_p)\, x_0} (\overline{a_0\ldots a_p}) \cdot a_0\ldots  
a_p$ sur la compl\'etion de Hilbert de $H_1^{dyn}$.   
Cela induit une repr\'esentation de la ${\rm C}^*$-alg\`ebre    
${\rm C}_0 (\mX_\Gamma,{\mathcal E})$ des sections du fibr\'e   
${\mathcal E}= ( {\rm C}(\Lambda_\Gamma)  \times \H^3)/\Gamma \longrightarrow  
\mX_\Gamma$, o\`u $\Gamma$ agit diagonalment sur  
${\rm C}(\Lambda_\Gamma) \times \H^3$.  
  
Nous consid\'erons l'op\'erateur lin\'eaire illimit\'e $D: {\mathcal  
H}^1_{dyn}\oplus {\mathcal H}_1^{dyn} \to {\mathcal   
H}^1_{dyn}\oplus {\mathcal H}_1^{dyn}$ qui agit comme une
multiplication par   
le poids $D: x\mapsto p\cdot x$, sur $gr^\Gamma_{2p} H^1_{dyn}$  
avec $p\leq 0$ et sur $gr^\Gamma_{2p} H_1^{dyn}$ avec $p\geq 1$.

\begin{thm}\label{dyn-SP3}  
Soit ${\mathcal H}={\mathcal H}^1_{dyn} \oplus {\mathcal M}  
\otimes_{{\rm C}_0 (\mX_\Gamma,{\mathcal E})} {\mathcal H}_1^{dyn}$,  
o\`u ${\mathcal M}$ est un bimodule r\'ealisant l'\'equivalence de Morita  
entre les ${\rm C}^*$-alg\`ebres ${\rm C}_0 (\mX_\Gamma,{\mathcal E})$  
et ${\mathcal A}:=( {\rm C}(\Lambda_\Gamma)\otimes {\rm C}_0(\H^3))
\rtimes \Gamma$.     
Alors, $({\mathcal A},{\mathcal H},\tilde D)$ est un triplet spectral  
avec $\tilde D|_{{\mathcal H}^1_{dyn}} = D$ et  
$\tilde D|_{{\mathcal M} \otimes{\mathcal H}_1^{dyn}} = 1\otimes  
D$.   
\end{thm}

\end{document}